\documentclass[a4paper,11pt]{article}
\usepackage[pdftex]{graphicx}
\usepackage{relsize}
\usepackage{amsmath}
\usepackage{amsfonts}
\usepackage{precalc1}
\usepackage{url}

\begin{document}
\newcommand{\defeq}{\mathrel{\mathop:}=}
\newcommand{\eqdef}{\mathrel{\mathop=}:}
\newcommand{\field}[1]{\mathbb{#1}}
\newcommand{\vect}[1]{\vec{\mathbf{#1}}}
\newcommand{\vectc}[1]{\vec{\mathbf{\mathcal{#1}}}}
\newcommand{\vectv}[2]{\vect{{#1}}_{#2}}
\newcommand{\scal}[2]{\vect{{#1}} \cdot \vect{{#2}}}
\newcommand{\ex}[1]{e^{#1}}
\newcommand{\dlt}[1]{\bigtriangleup #1}
\newcommand{\dm}[1]{\det\left(\textbf{#1}\right)}
\newcommand{\ud}[1]{\,\mathrm{d}{#1}}
\newcommand{\udm}[2]{\,\mathrm{d}_{#1}{#2}}
\newcommand{\vd}[2]{\,\mathrm{d}^{#1}\vect{{#2}}}
\newcommand{\vdm}[3]{\,\mathrm{d}_{#2}^{#1}\vect{{#3}}}
\newcommand{\four}[1]{\mathcal{F}[{#1}]}
\newcommand{\ifour}[1]{\mathcal{F}^{-1}[{#1}]}
\newcommand{\fouru}[1]{\mathcal{F}_1[{#1}]}
\newcommand{\bra}[1]{\langle #1|}
\newcommand{\ket}[1]{|#1\rangle}
\newcommand{\braket}[2]{\langle #1|#2\rangle}
\newcommand{\ketbra}[2]{|#1\rangle\langle#2|}
\newcommand{\brakop}[3]{\langle #1|#2|#3\rangle}
\newtheorem{thm}{Theorem}[section]
\newtheorem{prop}[thm]{Proposition}
\title{Notes on $SU(2)$ invariant \\
  absolutely continuous probability measures}
\author{Giuseppe Vitillaro \\
  \date{July 4, 2025\\v0.32}}
\maketitle
\begin{abstract}
  A measure independence property of Lebesgue
  measurable \emph{convex cones} of $\field{C}^2$,
  for $SU(2)$ transformations invariant continuous
  probability joint distributions over $\field{C}^2$,
  will be proved using the existence of the Haar probability
  measure for compact Hausdorff topological groups.
\end{abstract}
\section{Introduction}
My attention and curiosity has been caught from
a statement I read on a paper about the emergence
of the Born Rule in theory of quantum observation,
recently published on arXiv by a physicist
\cite[p.~10]{quantum}:
\begin{quote}
  We do not have any information about the state of the
  photon. That implies that the statistical distribution of $\alpha$
  and $\beta$ does not depend on the choice of a basis, or in
  other words, the distribution of $(\alpha,\beta)$ must be invariant
  under $SU(2)$ transforms.
\end{quote}

I had never thought to properties of $SU(2)$ transformations
invariant, continuous probability joint distributions over
$\field{C}^2$ and I was surprised to read the author assumed,
in the following lines, the probabilities of the Lebesgue
measurable subsets of $\field{C}^2$
\begin{equation*}
  \left\{
    (\alpha,\beta) \in \field{C}^2
    \suchthat
    |a||\alpha| > |b||\beta|
  \right\}
  \mbox{,\,\,\,\,}
  a,b \in \field{C}^2 \setminus \{(0,0)\}
\end{equation*}
are independent from the choice
of one of such $SU(2)$ invariant distributions
and can be evaluated, to obtain the Born Rule
\begin{equation*}
  P\Big[ |a||\alpha| > |b||\beta| \Big] =
    \frac{|a|^2}{|a|^2+|b|^2}
  \mbox{,\,\,\,\,}
  a,b \in \field{C}^2 \setminus \{(0,0)\}
\end{equation*}
using a gaussian shaped probability joint distribution in
$\field{C}^2$ \cite[p.~10, eq.~(65)/(66)]{quantum},
that looks like just as a convenient choice to make simple
the integral evaluation.

Quoting the author words
\begin{quote}
  One possible distribution that realizes this symmetry
  is given by
  \begin{eqnarray*}
    && \alpha \defeq G_1 + iG_2\\
    && \beta  \defeq G_3 + iG_4
  \end{eqnarray*}
  Where $G_n$ are mutually independent identically distributed
  gaussian random variables with zero mean. Because
  the sum of two independent Gaussian variables
  results in a new Gaussian variable, this construction
  is invariant under unitary transformations.
\end{quote}

In these notes I will just try to check if this 
claim about $SU(2)$ transformations invariant
continuous probability joint distributions is correct,
justifying the procedure the author has been following
to evaluate the Born Rule with an arbitrary, convenient,
choice of a probability distribution over $\field{C}^2$.

\section{Coordinates transformations}
We are going to work over the ``cartesian coordinates domain'' of couples
$\vect{w}=(\vect{\alpha},\vect{\beta})=(x+iy,u+iv)=(x,y,u,v) \in
\field{C}^2=\field{R}^4$
\begin{equation*}
  D_c=\left\{
    (\vect{\alpha},\vect{\beta}) \in \field{C}^2 
    \suchthat
    \vect{\alpha} \ne \vect{0} \text{\,\,or\,\,} \vect{\beta} \ne \vect{0}    
  \right\}=
  \C^2 \setminus \left\{(\vect{0},\vect{0})\right\}
\end{equation*}
for which we will introduce some suitable coordinates transformations.\\

The differentiable map $T_d:D_d \to D_c$
\begin{equation}
    \label{eq:dpolar}
    T_d(r,\rho,\phi,\theta) = (r\cos\phi,r\sin\phi,\rho\cos\theta,\rho\sin\theta)
\end{equation}
over the ``double polar coordinates domain''
\begin{equation*}
  D_d=\left\{
    (r,\rho,\phi,\theta) \in \field{R}^4
    \suchthat
    r \ge 0, \rho \ge 0, r>0\text{\,\,or\,\,}\rho>0\text{,\,\,}
    0 \le \phi,\theta \le 2\pi
  \right\}
\end{equation*}
has, almost everywhere, inverse map $T_d^{-1}:D_c \to D_d$
\begin{equation*}
  T_d^{-1}(\vect{\alpha},\vect{\beta})=
  T_d^{-1}(x,y,u,v) = (\sqrt{x^2+y^2},\sqrt{u^2+v^2},
    \eta(y,x),\eta(v,u))
\end{equation*}
where $\eta$ is the bijective $[0,2\pi)$ version of the atan2 function, 
as defined in numerical libraries.\\
 
The coordinates transformation $T_d$ has
an almost everywhere positive defined absolute value Jacobian
determinant
\begin{equation*}
  \begin{array}{ll}
    |\det(J_{T_d})|=r\rho\\
    \ud{x}\ud{y}\ud{u}\ud{v} = 
    |\det(J_{T_d})|\ud{r}\ud{\rho}\ud{\phi}\ud{\theta} =
    r\rho\ud{r}\ud{\rho}\ud{\phi}\ud{\theta}
  \end{array}
\end{equation*}
being $\vd{4}{w}=\vd{2}{\alpha}\vd{2}{\beta}=\ud{x}\ud{y}\ud{u}\ud{v}$ and 
$\ud{r}\ud{\rho}\ud{\phi}\ud{\theta}$
the hypervolume elements of $D_c$ and $D_d$,
in the sense of the Lebesgue integral.\\

Intuition suggests now to use a standard polar coordinates $(l,\psi)$
transformation of the radial coordinates
$(r,\rho)\text{,\,\,\,}r>0\text{\,\,or\,\,}\rho>0$
\begin{equation*}
  \left\{ 
    \begin{array}{ll}
      r = l\cos\psi\\ 
      \rho = l\sin\psi
    \end{array}
   \right. 
   \left.  
     \begin{array}{ll}
       l > 0 \\
       0 \le \psi \le \frac{\pi}{2}
     \end{array}
   \right. 
\end{equation*}\\
with $l^2=r^2+\rho^2>0$, $0 \le \psi \le \frac{\pi}{2}$,
$\psi=\arctan\frac{\rho}{r}$
\footnote
{
  Where it will appear $\arctan(r/s)$, with $r,s$
  non negative real numbers and $r>0$, we will
  take the freedom to think $\arctan(r/0)=\frac{\pi}{2}$,
  for the sake of the notation.
},
$\ud{r}\ud{\rho}=l\ud{l}\ud{\psi}$
to define a modified version of ``hyperspherical coordinates'',
suitable for our goals.\\

The differentiable map $T_h:D_h \to D_c$
\begin{equation}
  \label{eq:hyper}
    T_h(r,\rho,\phi,\theta) =
    (
      l\cos\psi\cos\phi,l\cos\psi\sin\phi,
      l\sin\psi\cos\theta,l\sin\psi\sin\theta
    )
\end{equation}
over the ``hyperspherical coordinates domain''
$(0,+\infty) \times [0,\frac{\pi}{2}] \times [0,2\pi]
\times [0,2\pi]$
\begin{equation*}
  D_h=\left\{
    (l,\psi,\phi,\theta) \in \field{R}^4
    :
    l > 0,
    0 \le \psi \le \frac{\pi}{2},
    0 \le \phi,\theta \le 2\pi
  \right\}
\end{equation*}
being the angular coordinates $(\psi,\phi,\theta)$
the ``Hopf coordinates'' \cite{hopf} imbedding $S^3$ in $\field{C}^2$,
has, almost everywhere, inverse map $T_h^{-1}:D_c \to D_h$
\begin{multline*}
  T_h^{-1}(\vect{\alpha},\vect{\beta})=
  T_h^{-1}(x,y,u,v) = \\
  = \left(
    \sqrt{x^2+y^2+u^2+v^2},
    \arctan{\sqrt{\frac{u^2+v^2}{x^2+y^2}}},
    \eta(y,x),\eta(v,u)
  \right)
\end{multline*}
 
Since $l>0$ and $0 \le \psi \le \frac{\pi}{2}$
over $D_h$, the coordinates transformation $T_h$ has
an almost everywhere positive defined absolute value Jacobian
determinant
\begin{equation*}
  \begin{array}{ll}
    |\det(J_{T_h})|=|-l^3\sin\psi\cos\psi|=
    l^3\sin\psi\cos\psi=\frac{1}{2}l^3\sin2\psi\\ 
    \ud{x}\ud{y}\ud{u}\ud{v} = 
    |\det(J_{T_h})|\ud{r}\ud{\rho}\ud{\phi}\ud{\theta} =
    \frac{1}{2}l^3\sin2\psi\ud{l}\ud{\psi}\ud{\phi}\ud{\theta}
  \end{array}
\end{equation*}
being $\vd{4}{w}=\vd{2}{\alpha}\vd{2}{\beta}=\ud{x}\ud{y}\ud{u}\ud{v}$ and 
$\ud{l}\ud{\psi}\ud{\phi}\ud{\theta}$
the hypervolume elements of $D_c$ and $D_h$,
again in the sense of the Lebesgue integral.\\

Note that in our notation $r=|\vect{\alpha}|$, $\rho=|\vect{\beta}|$ and
$l=\sqrt{|\vect{\alpha}|^2+|\vect{\beta}|^2}$, thinking to
$(\vect{\alpha},\vect{\beta})$ as a vector in $\field{C}^2$.

\section{Integrals evaluation rules}

We will use these coordinates tranformations to evaluate
Lebesgue integrals between their coordinates domains, beside
a Lebesgue measurable \emph{null set}, using them for
a change of variables and applying a well known theorem
\cite[p.~173, Theorem~8.26]{rudin},
\cite[p.~300, Theorem~263D]{fremlin2}.

If $p_d$ is a real valued, non negative, integrable function over
the double polar coordinates domain $D_d$ then
\begin{equation*}
  p_c(\vect{\alpha},\vect{\beta}) = p_c(x,y,u,v) = p_d({T_d}^{-1}(x,y,u,v)))
\end{equation*}
is a real valued, non negative, integrable function over the
cartesian coordinates domain $D_c$
and for every measurable subset $A$ of $D_c$
\begin{equation}
  \label{eq:ruled}
  \iiiint\limits_{A}p_c(\vect{\alpha},\vect{\beta})\vd{2}{\alpha}\vd{2}{\beta} =
  \iiiint\limits_{{T_d}^{-1}(A)} r\rho p_d(r,\rho,\phi,\theta) 
     \ud{r}\ud{\rho}\ud{\phi}\ud{\theta}
\end{equation}\\
From the other side if $p_c(\alpha,\beta)$ is a real valued
non negative,integrable
function over the cartesian coordinates domain $D_c$ then
\begin{equation*}
  p_d(r,\rho,\phi,\theta) = p_c(T_d(x,y,u,v))
\end{equation*}
is a real valued, non negative integrable function over the
double polar coordinates domain $D_d$
and for every measurable subset $A$ of $D_d$
\begin{equation*}
  \iiiint\limits_{A}p_d(r,\rho,\phi,\theta)\ud{r}\ud{\rho}\ud{\phi}\ud{\theta} =
  \iiiint\limits_{{T_d}(A)}\frac{1}{|\vect{\alpha}||\vect{\beta}|}
    p_c(\vect{\alpha},\vect{\beta})\vd{2}{\alpha}\vd{2}{\beta}
\end{equation*}

If $p_h$ is a real valued, non negative, integrable function
over the hyperspherical polar coordinates
domain $D_h$ then
\begin{equation*}
  p_c(\vect{\alpha},\vect{\beta}) = p_c(x,y,u,v) = p_h({T_h}^{-1}(x,y,u,v)))
\end{equation*}
is a real valued, non negative, integrable function over
the cartesian coordinates domain $D_c$
and for every measurable subset $A$ of $D_c$
\begin{equation}
  \label{eq:ruleh}
  \iiiint\limits_{A}p_c(\vect{\alpha},\vect{\beta})\vd{2}{\alpha}\vd{2}{\beta} =
  \frac{1}{2}\iiiint\limits_{{T_h}^{-1}(A)} l^3 p_h(l,\psi,\phi,\theta)
     \sin2\psi
     \ud{l}\ud{\psi}\ud{\phi}\ud{\theta}
\end{equation}\\
Again if $p_c(\alpha,\beta)$ is a real valued, non negative,
integrable function over the cartesian coordinates domain $D_c$ then
\begin{equation*}
  p_h(l,\psi,\phi,\theta) = p_c(T_h(x,y,u,v))
\end{equation*}
is a real valued, non negative, integrable function over
the double polar coordinates domain $D_h$
and for every measurable subset $A$ of $D_h$
\begin{equation*}
  \iiiint\limits_{A}p_h(l,\psi,\phi,\theta)\ud{l}\ud{\psi}\ud{\phi}\ud{\theta} =
  \iiiint\limits_{{T_h}(A)}\frac{1}{
    |\vect{\alpha}||\vect{\beta}|
    \sqrt{|\vect{\alpha}|^2 + |\vect{\beta}|^2 }
  }
  p_c(\vect{\alpha},\vect{\beta})\vd{2}{\alpha}\vd{2}{\beta}
\end{equation*}

We may want to use these coordinates transformations and rules to
evaluate integrals over the coordinates domains we have defined.

\section{A list of simple integrals}
Some simple trigonometric and exponential integral
evaluation we will find useful in the following sections
\begin{alignat}{4}
  \label{eq:trig1}
    &\int_{0}^{x}\sin2\psi\ud{\psi}=\sin^2x
    \text{,\,\,\,\,}
    0 \le x \le \frac{\pi}{2}\\
  \label{eq:trig2}
    &\int_{0}^{\frac{\pi}{2}}\sin2\psi\ud{\psi}=1\\
  \label{eq:exp1}
    &\int_{0}^{l}q\ex{-\frac{q}{2}}\ud{q}=
    4-2(l+2)\ex{-\frac{q}{2}}
    \text{,\,\,\,\,}
    l \ge 0\\
  \label{eq:exp2}
    &\int_{0}^{+\infty}q\ex{-\frac{q}{2}}\ud{q}=4
\end{alignat}

\section{$SU(2)$ invariant functions}

The \emph{special unitary} group $SU(2)$ of degree 2 is defined as the
group, under the operation of matrix multiplication, of
$2 \times 2$ complex unitary matrices with determinant equal to 1:
\begin{equation*}
  SU(2) = \left\{
           \begin{pmatrix}
	     \vect{\alpha} &  -\overline{\vect{\beta}}\\
	     \vect{\beta}  &   \overline{\vect{\alpha}}
	   \end{pmatrix}
           \text{\,}\Bigg|\text{\,\,\,\,}
	   \vect{\alpha},\vect{\beta} \in \field{C},
	   \text{\,\,\,\,}
           \vect{\alpha}\overline{\vect{\alpha}} +
	   \vect{\beta}\overline{\vect{\beta}} = 1
	  \right\}
\end{equation*}
Each matrix $U \in SU(2)$ defines, under the operation of matrix
vector multiplication, an orthogonal transformation
$T_U:\field{C}^2 \to \field{C}^2$, $T_U(\vect{w})=U\vect{w}\text{,\,\,\,\,}
\vect{w} \in \field{C}^2$:
\begin{equation*}
  \langle T_U(\vect{w}_1),T_U(\vect{w}_2) \rangle =
  \langle U\vect{w}_1,U\vect{w}_2 \rangle =
  \langle \vect{w}_1,\vect{w}_2 \rangle
  \text{,\,\,\,\,}
  \vect{w}_1,\vect{w}_2 \in \field{C}^2
\end{equation*}
such that $|\vect{w}|=|T_U(\vect{w})|=|U\vect{w}|$, for each
$\vect{w} \in \field{C}^2$.\\

Suppose $p_c:\field{C}^2 \to [0,+\infty)$ is a real valued,
non negative function over $\field{C}^2$ for which there exists
a real valued, non negative, 
function $f:[0,+\infty) \to [0,+\infty)$ such that
$p_c(\vect{w})=f(|\vect{w}|)$, for each $\vect{w} \in \field{C}^2$.\\

Then, for each $U \in SU(2)$:
\begin{equation*}
  p_c(U\vect{w})=f(|U\vect{w}|)=f(|\vect{w}|)=p_c(\vect{w})
  \text{,\,\,\,\,}
  \vect{w} \in \field{C}^2
\end{equation*}
and hence $p_c$ is invariant for $SU(2)$ transformations,
i.e. $p_c(U\vect{w})=p_c(\vect{w})$, for each $\vect{w} \in \field{C}^2$
and $U \in SU(2)$.\\

Now suppose the real valued, non negative function
$p_c:\field{C}^2 \to \field{R}$
over $\field{C}^2$ is invariant for $SU(2)$ transformations.\\

For each $\vect{w}=(\vect{\alpha},\vect{\beta}) \in \field{C}^2$,
with $\vect{w} \ne (\vect{0},\vect{0})$ we may define the $2 \times 2$
complex matrix:
\begin{equation}
  \label{eq:rotmat}
  U_{\vect{w}}=U_{(\vect{\alpha},\vect{\beta})}=
  \begin{pmatrix}
    \overline{\vect{\alpha}}/l & \overline{\vect{\beta}}/l \\
    \text{\,\,\,\,} & \text{\,\,\,\,} \\
    -\vect{\beta}/l & \vect{\alpha}/l
  \end{pmatrix}
\end{equation}
where $l^2=\vect{\alpha}\overline{\vect{\alpha}} +
\vect{\beta}\overline{\vect{\beta}}$ and $l=|\vect{w}| > 0$.\\

It is easy to verify $U_{\vect{w}}$ is an $SU(2)$ matrix
such that
\begin{equation*}
  U_{\vect{w}}\vect{w}=U_{(\vect{\alpha},\vect{\beta})}
  (\vect{\alpha},\vect{\beta})=(l,\vect{0})=(|\vect{w}|,\vect{0})
\end{equation*}

Using this unitary matrices family, we may now verify the real valued,
non negative function, over the non negative real
numbers $f:[0,+\infty) \to [0,+\infty)$ defined by
\begin{equation*}
  f(l) = p_c(l,\vect{0}) \mbox{,\,\,\,\,} l \ge 0
\end{equation*}
satisfy the properties $p_c(\vect{w})=p_c(U_{\vect{w}}\vect{w})=f(|\vect{w}|)$,
for each $\vect{w} \in \field{C}^2$, $\vect{w} \ne (\vect{0},\vect{0})$ and
$p_c(\vect{0},\vect{0})=f(0)$, just because the function $p_c$ is
invariant for this particular unitary transformation $U_{\vect{w}}$
by assumption.\\

We have proved:
\begin{prop}
  \label{inv1}
  A real valued, non negative function $p_c:\field{C}^2 \to [0,+\infty)$
  is invariant for $SU(2)$ transformations \emph{if and only if} there
  exists a function $f:[0,+\infty) \to [0,+\infty)$ such
  that $p_c(\vect{w})=f(|\vect{w}|)$, for each
  $\vect{w} \in \field{C}^2$.\\
\end{prop}

\section{$SU(2)$ invariant integrable functions}

From now we will use $\mathcal{L}(\field{C}^2)$ to denote
the Lebesgue complete \emph{sigma algebra} of Lebesgue measurable subsets
of $\field{C}^2$ and $\lambda$ for the Lebesgue measure over
$\mathcal{L}(\field{C}^2)$. We will think to integrals of
functions as Lebesgue integrals and we will denote
$\vd{4}{w}=\vd{2}{\alpha}\vd{2}{\beta}$
with $\vdm{4}{\lambda}{w}$ or $\ud{\lambda}$.\\

A property will hold \emph{almost everywhere} or a.e.
over $\field{C}^2$, if the set of elements where the property
does not hold is a Lebesgue ``\emph{null set}'', i.e. it is
a Lebesgue measurable subset of $\field{C}^2$ and its
Lebesgue measure is \emph{zero}.\\

A Lebesgue measurable function $p_c:\field{C}^2 \to [0,+\infty)$
will be called ``\emph{invariant a.e. for $SU(2)$ transformations}''
if for every $U \in SU(2)$ there exists a Lebesgue null set $N_U$,
possibly depending from $U$, such that $p_c(U\vect{w})=p_c(\vect{w})$,
for each $\vect{w} \in \field{C}^2$, $\vect{w} \notin N_U$, where
the family of $\field{C}^2$ null sets $\left\{N_U\right\}_{U \in SU(2)}$
is possibly uncountable.\\

Let $p_c:\field{C}^2 \to [0,+\infty)$ a measurable function for
which there exists a real valued, non negative, measurable function
$f:[0,+\infty) \to [0,+\infty)$ such that
$p_c(\vect{w})=f(|\vect{w}|)$ a.e. over $\field{C}^2$,
i.e. there exists a Lebesgue null set $N$ such that
$p_c(\vect{w})=f(|\vect{w}|)$ for each $\vect{w} \in \field{C}^2,
\vect{w} \notin N$.

For any matrix $U \in SU(2)$ and for any $\vect{w} \notin N$,
being $|U\vect{w}|=|\vect{w}|$, it holds
$p_c(U\vect{w})=f(|U\vect{w}|)=f(|\vect{w}|)=p_c(\vect{w})$,
hence $p_c$ is invariant a.e. for $SU(2)$ transformations.\\

Now let suppose $p_c:\field{C}^2 \to [0,+\infty)$
is a Lebsgue measurable, integrable function, invariant a.e. for $SU(2)$
transformations.\\

Thinking to a well know diffeomorphism we will denote
$SU(2)$ with the same symbol $S^3$ used for the
3-sphere, just as a useful shortcut in notation.\\

Because $SU(2)$ is a compact Hausdorff topological group
\cite[p.~280, 441X(h)(ii)]{fremlin4i} it has an
\emph{Haar probability} $\nu$, i.e. a \emph{totally finite}
measure over the sigma algebra generated by its topology,
invariant for $SU(2)$ transformations \cite{haarw}:
$\nu$ is a \emph{normalized} positive measure, $\nu(S^3)=1$, and
$\nu(UA)=\nu(\{U\vect{w} \suchthat \vect{w} \in A\})=\nu(A)$,
for each measurable subset $A$ of $SU(2)$.\\

We will denote Haar integrals of real valued, 
measurable functions over $S^3$ with the symbols
$\udm{\nu}{U}$ or $\ud{\nu}$ and observe for
every real valued, $\nu$-integrable function
$f:S^3 \mapsto \field{R}$ holds
\begin{equation}
  \label{eq:haarinv}
  \int\limits_{S^3}
  f(VU)\udm{\nu}{U}=
  \int\limits_{S^3}
  f(U)\udm{\nu}{U}
  \mbox{,\,\,\,\,for each }
  V \in S^3
\end{equation}

Because the Haar probability $\nu$ is finite and the
Lebsgue measure $\lambda$ is $\sigma$-finite, there exists
one and only one \emph{product measure}, $\udm{\nu}{U}\vdm{4}{\lambda}{w}$,
over the product space $S^3 \times \field{C}^2$
\cite[p.~114, Theorem~10.4]{bartle2} and the
the map $(U,\vect{w}) \mapsto U\vect{w}$,
$(U,\vect{w}) \in S^3 \times \field{C}^2$ is
$\udm{\nu}{U}\vdm{4}{\lambda}{w}$-measurable.

Then for each real valued, non negative, measurable,
bounded function $h:\field{C}^2 \to [0,+\infty)$
it makes sense to look to the integral
\begin{equation*}
  \iint\limits_{S^3 \times \field{C}^2}
    p_c(U\vect{w})h(\vect{w})
    \udm{\nu}{U}\vdm{4}{\lambda}{w}
\end{equation*}

For the Tonelli's Theorem \cite[p~.118, Theorem~10.9]{bartle2},
being $(U,\vect{w}) \mapsto p_c(U\vect{w})h(\vect{w})$,
$(U,\vect{w}) \in S^3 \times \field{C}^2$ a real valued, non negative,
measurable function over $S^3 \times \field{C}^2$ and $S^3$,
$\field{C}^2$ $\sigma$-finite measure spaces
\begin{equation}\begin{gathered}
  \label{eq:tonelli}
  \iint\limits_{S^3 \times \field{C}^2}
    p_c(U\vect{w})h(\vect{w})
    \udm{\nu}{U}\vdm{4}{\lambda}{w} =\\
  =
  \int\limits_{S^3}
  \left[
    \int\limits_{\field{C}^2}
      p_c(U\vect{w})h(\vect{w})
      \vdm{4}{\lambda}{w}
  \right]
      \udm{\nu}{U}=
  \int\limits_{\field{C}^2}
  \left[
    \int\limits_{S^3}
      p_c(U\vect{w})
      h(\vect{w})
     \udm{\nu}{U}
  \right]
    \vdm{4}{\lambda}{w}
\end{gathered}\end{equation}
where all functions under integration are non negative,
measurable functions, with values in $[0,+\infty]$,
possibly assuming the non real value $+\infty$.

But we know that, by assumption, for each $U \in S^3$ the function
$p_c(U\vect{w})$ is equal a.e. to the function $p_c(\vect{w})$
over $\field{C}^2$ and we may use this property to
evaluate the second partial integral in $\eqref{eq:tonelli}$
\begin{equation*}
  \int\limits_{\field{C}^2}
    p_c(U\vect{w})h(\vect{w})
    \vdm{4}{\lambda}{w}=
  \int\limits_{\field{C}^2}
    p_c(\vect{w})h(\vect{w})
    \vdm{4}{\lambda}{w}
\end{equation*}
for every $U \in S^3$ and because $p_c$ is a real valued, non negative,
integrable function and $h$ a real valued, non negative, measurable, bounded
function, over their domain $\field{C}^2$, the function
$\vect{w} \mapsto p_c(\vect{w})h(\vect{w})$, $\vect{w} \in \field{C}^2$
is integrable over $\field{C}^2$ \cite[p~.48, Exercise~5.E]{bartle2}
and hence the function 
\begin{equation*}
    U \mapsto
    \int\limits_{\field{C}^2}
    p_c(U\vect{w})h(\vect{w})
    \vdm{4}{\lambda}{w}=
    \int\limits_{\field{C}^2}
    p_c(\vect{w})h(\vect{w})
    \vdm{4}{\lambda}{w} \in \field{R}
    \mbox{,\,\,\,\,}U \in S^3
\end{equation*}
is a constant, non negative, real valued function, $\nu$-integrable
over $S^3$, being $\nu$ a finite measure.

This imply the function $(U,\vect{w}) \mapsto p_c(U\vect{w})h(\vect{w})$,
$(U,\vect{w}) \in S^3 \times \field{C}^2$,
is $\udm{\nu}{U}\vdm{4}{\lambda}{w}$ integrable,
globally over the product space $S^3 \times \field{C}^2$,
and its integral is
\begin{multline}
  \label{eq:inteq}
   \iint\limits_{S^3 \times \field{C}^2}
    p_c(U\vect{w})h(\vect{w})
    \udm{\nu}{U}\vdm{4}{\lambda}{w}=
   \int\limits_{S^3}
   \left[
     \int\limits_{\field{C}^2}
     p_c(\vect{w})h(\vect{w})
     \vdm{4}{\lambda}{w}
   \right]
      \udm{\nu}{U}=\\
   =
   \nu(S^3)
     \int\limits_{\field{C}^2}
     p_c(\vect{w})h(\vect{w})
     \vdm{4}{\lambda}{w}=
   \int\limits_{\field{C}^2}
     p_c(\vect{w})h(\vect{w})
     \vdm{4}{\lambda}{w} < +\infty
\end{multline}
being $\nu$ a normalized
probabilty measure over $S^3$ with $\nu(S^3)=1$,
for each real valued, non negative, measurable, bounded
function $h:\field{C}^2 \to [0,+\infty)$.\\

We are now in the condition to apply the Fubini's Theorem
\cite[p.~119, Theorem~10.10]{bartle2}
to the last partial integral in $\eqref{eq:tonelli}$,
with the constant function $h:\vect{w} \mapsto 1,
\vect{w} \in \field{C}^2$, to reach the
conclusion the \emph{extended real} valued, non negative,
measurable function over $\field{C}^2$
\begin{equation*}
  \vect{w} \mapsto
  \int\limits_{S^3}
     p_c(U\vect{w})
     \udm{\nu}{U}
  \mbox{,\,\,\,\,}\vect{w} \in \field{C}^2
\end{equation*}
is equal \emph{almost everywhere} to a \emph{real} valued, non negative,
\emph{integrable} function over $\field{C}^2$.

In other words there exists a \emph{null set}
$N \in \mathcal{L}(\field{C}^2)$, \emph{independent}
from $SU(2)$ unitary transformations, such that
$U \mapsto p_c(U\vect{w})$, $U \in S^3$, 
is $\nu$-integrable over $S^3$
\begin{equation*}
  0 \le 
  \int\limits_{S^3}
     p_c(U\vect{w})
     \udm{\nu}{U} < +\infty
\end{equation*}
for every $\vect{w} \in \field{C}^2$, $\vect{w} \notin N$,
so the real valued, non negative, function
$g:\field{C}^2 \to [0,+\infty)$ is well defined by
\begin{equation*}
  g(\vect{w})=
  \begin{cases}
     \mathlarger{\int\limits_{S^3}}
     p_c(U\vect{w})
     \udm{\nu}{U}
     &\mbox{if }\vect{w} \in \field{C}^2
     \mbox{,\,\,\,\,}\vect{w} \notin N\\[1cm]
     0  
     &\mbox{if }\vect{w} \in N
  \end{cases}
\end{equation*}
and is integrable over $\field{C}^2$
\begin{equation*}
  \int\limits_{\field{C}^2}
     g(\vect{w})\vdm{4}{\lambda}{w}=
  \int\limits_{\field{C}^2}
  \left[
    \int\limits_{S^3}
       p_c(U\vect{w})
       \udm{\nu}{U}
  \right]
  \vdm{4}{\lambda}{w}=
  \int\limits_{\field{C}^2}
    p_c(\vect{w})
    \vdm{4}{\lambda}{w} < +\infty
\end{equation*}

The function $g$, using again $\eqref{eq:tonelli}$,
for each unitary matrix $V \in SU(2)$, has the nice property
\begin{equation*}
   g(\vect{w})=
   \int\limits_{S^3}
     p_c(U\vect{w})
     \udm{\nu}{U}=
   \int\limits_{S^3}
     p_c(VU\vect{w})
     \udm{\nu}{U}=
   g(V\vect{w})
\end{equation*}
for every $\vect{w} \in \field{C}^2$, $\vect{w} \notin N$,
for the property $\eqref{eq:haarinv}$ of real valued
$\nu$-integrable functions, being the map
$U \mapsto p_c(U\vect{w})$, $U \in S^3$,
$\nu$-integrable over $S^3$.

In particular this holds for the unitary matrix $U_{\vect{w}}$
we have define in $\eqref{eq:rotmat}$ such that
$g(U_{\vect{w}}\vect{w})=g(|\vect{w}|,\vect{0})$,
for every $\vect{w} \in \field{C}^2$, $\vect{w} \notin N$.

Then the real valued, non negative, function
$f:[0,+\infty) \to [0,+\infty)$ defined by
$f:l \mapsto g(l,\vect{0})$, $l \ge 0$ is measurable and
\begin{equation*}
  g(\vect{w})=g(U_{\vect{w}}\vect{w})=
  g(|\vect{w}|,\vect{0})=f(|\vect{w}|)
\end{equation*}
for every $\vect{w} \in \field{C}^2$, $\vect{w} \notin N$,
i.e. $g(\vect{w})=f(|\vect{w}|)$ a.e. over $\field{C}^2$.

We may now combine the definition of the function $g$
and the equations $\eqref{eq:tonelli}$,
$\eqref{eq:inteq}$ to get
\begin{equation*}
  \int\limits_{\field{C}^2}
     p_c(\vect{w})h(\vect{w})
     \vdm{4}{\lambda}{w}
   =
 \int\limits_{\field{C}^2}
   h(\vect{w})
   \left[
    \int\limits_{S^3}
      p_c(U\vect{w})
     \udm{\nu}{U}
   \right]
   \vdm{4}{\lambda}{w}
  =
  \int\limits_{\field{C}^2}
    g(\vect{w})h(\vect{w})
    \vdm{4}{\lambda}{w}
\end{equation*}
for every real valued, non negative, measurable, bounded
function $h$ over $\field{C}^2$ and this imply,
choosing suitable measurable, non negative, bounded functions
$h^+$, $h^-$, defined in terms of the Lebesgue integrable function
$[p_c(\vect{w})-g(\vect{w})]$, and
applying a simple measure theory property
\cite[p.~34, Corollary~4.10]{bartle2}, that
$p_c(\vect{w})=g(\vect{w})$
\emph{almost everywhere} over $\field{C}^2$ and
from this we obtain  $p_c(\vect{w})=f(|\vect{w}|)$
a.e. over $\field{C}^2$.\\

We have proved
\begin{prop}
  \label{inv2}
  A real valued, non negative, integrable
  function $p_c:\field{C}^2 \to [0,+\infty)$
  is invariant for $SU(2)$ transformations a.e.
  \emph{if and only if} there
  exists a measurable function $f:[0,+\infty) \to [0,+\infty)$ such
  that $p_c(\vect{w})=f(|\vect{w}|)$, a.e. over $\field{C}^2$.
\end{prop}
In other words the Proposition \ref{inv1} still holds
if \emph{everywhere} is replaced with \emph{almost everywhere}
over $\field{C}^2$ and we observe one of the implications
is strictly connected to the existence of an Haar probability
measure over the compact Hausdorff topological group $SU(2)$.\\

It is interesting to quote here a comment of Prof. D.H.Fremlin
\cite[p.~282, Notes 441]{fremlin4i}
about the \emph{general proof} of the existence of the Haar Measure
employs the \emph{Axiom of Choice}, but
``\emph{it can be built with much weaker principles}'',
i.e. the \emph{Axiom of Countable Choice} \cite{aocc},
see \cite[\S~27.4]{naimark}.\\

The sketch of the proof of Proposition \ref{inv2} is
essentially due to Prof. D.H. Fremlin
(personal communication, May 19, 2012).\\

\section{$SU(2)$ invariant measures}

Now we choose a finite, normalizable, 
measure $\mu$ over the Lebesgue sigma algebra
$\mu:\mathcal{L}(\field{C}^2) \to [0,+\infty)$,
$\mu(\field{C}^2) < +\infty$,
\emph{absolutely continuous} with respect the Lebesgue measure
$\lambda$, $\mu \ll \lambda$.

From the Radon-Nikod\'ym Theorem \cite[p.~85, Theorem~8.9]{bartle2}
we know there exists a real valued, non negative, integrable
function $p_c:\field{C}^2 \to [0,+\infty)$ such that
\begin{equation}
  \label{eq:murad}
  \mu(A)=\int\limits_{A}p_c(\vect{w})\vdm{4}{\lambda}{w}
\end{equation}
for each measurable set $A \in \mathcal{L}(\field{C}^2)$ and
that such a function $p_c$ is uniquely determined
$\lambda$-almost everywhere over $\field{C}^2$.
Its a.e. equivalence class is usually called the
\emph{Radon-Nikod\'ym derivative} $d\mu/d\lambda$ of the
measure $\mu$ with respect to $\lambda$.\\

For each unitary matrix $U \in SU(2)$ the corresponding
orthogonal transformation
$T_U:\field{C}^2 \to \field{C}^2$,
$T_U(\vect{w})=U\vect{w}\text{,\,\,\,\,}\vect{w} \in \field{C}^2$
is a diffeomorphism and, by definition, $|\det(J_{T_U})|=1$
over $\field{C}^2$, hence we may use $T_U$ as a coordinates
transformation map for a change of variables in the integral
\begin{equation}\begin{split}
  \label{eq:muinv}
  \mu(UA)&=
  \int\limits_{UA}p_c(\vect{w})\vdm{4}{\lambda}{w}=
  \int\limits_{T_U^{-1}(UA)}p_c(T_U(\vect{w}))|\det(J_{T_U})|
    \vdm{4}{\lambda}{w}=\\
  &=
  \int\limits_{A}p_c(U\vect{w})\vdm{4}{\lambda}{w}
\end{split}\end{equation}
to evaluate $\mu(UA)$ for each measurable set
$A \in \mathcal{L}(\field{C}^2)$ as a Lebesgue integral
of the function $p_c(U\vect{w})$.

Let suppose now the measure $\mu$ is invariant for SU(2)
transformations, i.e. $\mu(UA)=\mu(A)$, for each
unitary matrix $U \in SU(2)$ and for each
Lebesgue measurable set $A \in \mathcal{L}(\field{C}^2)$.

We can combine equations $\eqref{eq:murad}$, $\eqref{eq:muinv}$
to get
\begin{equation*}
  \int\limits_{A}p_c(\vect{w})\vdm{4}{\lambda}{w}=
  \mu(A)=
  \mu(UA)=
  \int\limits_{A}p_c(U\vect{w})\vdm{4}{\lambda}{w}
\end{equation*}
for each measurable set $A \in \mathcal{L}(\field{C}^2)$
and for each unitary matrix $U \in SU(2)$.

Then for each $U \in SU(2)$, the functions $\vect{w} \mapsto p_c(\vect{w})$
and $\vect{w} \mapsto p_c(U\vect{w})$, $\vect{w} \in \field{C}^2$,
have both to be Radon-Nikod\'ym
derivatives of the measure $\mu$ with respect to $\lambda$
and so, for the uniqueness a.e. of $d\mu/d\lambda$,
they have to be equal $p_c(\vect{w})=p_c(U\vect{w})$
almost everywhere, with respect the Lebesgue measure $\lambda$
over $\field{C}^2$.

In other words the Lebesgue integrable
function $p_c:\field{C}^2 \to [0,+\infty)$
is invariant for $SU(2)$ transformations a.e. and, for
Proposition \ref{inv2}, there exists a real valued,
non negative, Lebesgue measurable function
$f:[0,+\infty) \to [0,+\infty)$
such that $p_c(\vect{w})=f(|\vect{w}|)$ a.e. over
$\field{C}^2$ and, recalling the integral evaluation rule
$\eqref{eq:ruleh}$
\begin{equation}\begin{split}
  \mu(A)&=
  \int\limits_{A}p_c(\vect{w})\vdm{4}{\lambda}{w}=
  \int\limits_{A}f(|\vect{w}|)\vdm{4}{\lambda}{w}=
  \iiiint\limits_{A}f(|\vect{w}|)\vd{4}{w}=\\
  &=
  \frac{1}{2}
  \iiiint\limits_{T_h^{-1}(A)}
    l^3 f(l) \sin{2\psi}
    \ud{l}\ud{\psi}\ud{\phi}\ud{\theta}
\end{split}\end{equation}
for each Lebesgue measurable set $A \in \mathcal{L}(\field{C}^2)$.\\

We may want to evaluate such an integral over $\field{C}^2$
to use the fact $\mu$ is a \emph{finite measure}
\begin{equation}\begin{split}
 \label{eq:nder}
 \mu(\field{C}^2) &=
 \frac{1}{2}
 \iiiint\limits_{T_h^{-1}(\field{C}^2)}
   l^3 f(l) \sin{2\psi}
   \ud{l}\ud{\psi}\ud{\phi}\ud{\theta}=\\
 &=
 \frac{1}{2}
 \int_{0}^{+\infty}\int_{0}^{\frac{\pi}{2}}\int_{0}^{2\pi}\int_{0}^{2\pi}
    l^3 f(l) \sin{2\psi}
    \ud{l}\ud{\psi}\ud{\phi}\ud{\theta}=\\
 &=
 \frac{1}{2}
 \left[
   \int_{0}^{\frac{\pi}{2}}
     \sin{2\psi}\ud{\psi}
 \right]
 \left[
   \int_{0}^{2\pi}\int_{0}^{2\pi}
     \ud{\phi}\ud{\theta}
 \right]
 \left[
    \int_{0}^{+\infty}
      l^3 f(l) \ud{l}
 \right]=\\
 &=
 \frac{1}{2}\Big[1\Big]\Big[4\pi^2\Big]
 \left[
   \int_{0}^{+\infty}l^3 f(l)\ud{l}
 \right]=
 2\pi^2\int_{0}^{+\infty}l^3 f(l)\ud{l}
\end{split}\end{equation}\\
applying the Tonelli's Theorem,
to measurable non negative functions over their domains,
and using the integral $\eqref{eq:trig2}$.

We have discovered the real valued, non negative, measurable
function $l^3f(l)$ has to be Lebesgue
integrable over $[0,+\infty)$
\begin{equation}
  M_3(f)=\int_{0}^{+\infty}l^3 f(l)\ud{l}=
  \frac{1}{2\pi^2}\mu(\field{C}^2) < +\infty
\end{equation}
and, if $\mu$ is a \emph{probability measure}, i.e.
$\mu(\field{C}^2)=1$, we are actually able to evaluate
such an integral
\begin{equation}
  \label{eq:nconstr}
  M_3(f)=\int_{0}^{+\infty}l^3 f(l)\ud{l}=
  \frac{1}{2\pi^2} 
\end{equation}
being the value $1/2\pi^2$ the \emph{normalization constraint}
that let an absolutely continuous finite measure $\mu$
to be an absolutely continuous probability measure with a
Radon-Nikod\'ym derivative, with respect the Lebesgue
measure $\lambda$, $d\mu/d\lambda$ equal almost everywhere to 
$p_c(\vect{w})=f(|\vect{w}|)$, over $\field{C}^2$.\\

We will denote with $\mathcal{M}_3$ the class of real valued,
non negative, measurable functions $f:[0,+\infty) \to [0,+\infty)$
with finite \emph{third moment} $M_3(f) < +\infty$
and with $\mathcal{M}_{3U}$ its subclass of functions
satisfying the normalization constraint $\eqref{eq:nconstr}$,
$M_3(f)=1/2\pi^2$.\\

Let we choose an element $f \in \mathcal{M}_3$ and
let we define the real valued, non negative, measurable
function $p_c:\field{C}^2 \to [0,+\infty) \mbox{,\,\,\,\,}
p_c(\vect{w})=f(|\vect{w}|) \mbox{, } \vect{w} \in \field{C}^2$.

Then, just repeating the same derivations in
$\eqref{eq:nder}$, the function
$\mu:\mathcal{L}(\field{C}^2) \to [0,+\infty)$
defined by
\begin{equation*}
  \mu(A)=
  \int\limits_{A}p_c(\vect{w})\vdm{4}{\lambda}{w}=
  \int\limits_{A}f(|\vect{w}|)\vdm{4}{\lambda}{w}
  \mbox{,\,\,\,\,}
  A \in \mathcal{L}(\field{C}^2)
\end{equation*}
is an absolutely continuous, $\mu \ll \lambda$,
finite measure over $\field{C}^2$ and its Radon-Nikod\'ym
derivate $d\mu/d\lambda$ is equal almost everywhere
to the function $p_c$ over $\field{C}^2$.

If the function $f$ belongs to the subclass $\mathcal{M}_{3U}$,
i.e. it satisfies the normalization constraint, then
$\mu$ is actually a probability measure over
Lebesgue measurable sets $\mathcal{L}(\field{C}^2)$.

Furthermore, for each unitary matrix $U \in SU(2)$,
by the same argument used in $\eqref{eq:muinv}$
\begin{equation*}
  \mu(UA)=
  \int\limits_{UA}f(|\vect{w}|)\vdm{4}{\lambda}{w}=
  \int\limits_{A}f(|U\vect{w}|)\vdm{4}{\lambda}{w}=
  \int\limits_{A}f(|\vect{w}|)\vdm{4}{\lambda}{w}=
  \mu(A)
\end{equation*}
for each Lebesgue measurable set $A \in \mathcal{L}(\field{C}^2)$
and hence the absolutely continuous finite measure $\mu$ is
invariant for $SU(2)$ transformations.\\

We have proved
\begin{prop}
  \label{inv3m}
  An absolutely continuous finite measure
  $\mu:\mathcal{L}(\field{C}^2) \to [0,+\infty)$,
  $\mu \ll \lambda$,
  is invariant for $SU(2)$ transformations
  \emph{if and only if}
  there exists a function $f \in \mathcal{M}_3$
  such that the Radon-Nikod\'ym derivative $d\mu/d\lambda$
  of $\mu$ with respect the Lebesgue measure $\lambda$ is
  equal almost everywhere to $f(|\vect{w}|)$ over $\field{C}^2$,
  i.e. $(d\mu/d\lambda)(\vect{w})=f(|\vect{w}|)$, a.e.
  over $\field{C}^2$. Furthermore $\mu$ is a probability measure
  \emph{if and only if}
  $f \in \mathcal{M}_{3U}$, i.e.
  the function $f$ satisfies the normalization constraint
  $M_3(f)=1/2\pi^2$.
\end{prop}

In other words there is a \emph{bijective} correspondence
between $SU(2)$ invariant, absolutely continuous, finite
measures over $\mathcal{L}(\field{C}^2)$ and equivalence
classes of \emph{equal a.e} real valued,
non negative, measurable functions $f:[0,+\infty) \to [0,+\infty)$
with finite third moment
$M_3(f)=\int_{0}^{+\infty}l^3 f(l)\ud{l}<+\infty$,
defined by
\begin{eqnarray}
  \label{eq:corrmu}
  && \mu_f(A)=
  \frac{1}{2}
  \iiiint\limits_{T_h^-1(A)}
    l^3 f(l) \sin2\psi
    \ud{l}\ud{\psi}\ud{\phi}\ud{\theta}
  \mbox{,\,\,\,\,}
  A \in \mathcal{L}(\field{C}^2)\\
  \label{eq:corrpc}
  && p_c[\mu](\vect{w})=\frac{d\mu}{d\lambda}(\vect{w})
  \mbox{,\,\,\,\,a.e. over } \field{C}^2\\
  \label{eq:corrf}
  && f_\mu(l) = p_c[\mu](l,\vect{0})=
  \frac{d\mu}{d\lambda}(l,\vect{0})
  \mbox{,\,\,\,\,a.e. over } [0,+\infty)
\end{eqnarray}

This bijective correspondence, and the derived notation
we will continue to use in the following sections,
express the intuition an $SU(2)$ invariant, absolutely
continuous finite measure is well defined \emph{if and only if}
we know its Radon-Nikod\'ym derivative,
with respect the Lebesgue measure $\lambda$,
a.e. along \emph{exactly}
one radius of $\field{C}^2$. 

The measure $\mu$ is a probability measure
if and only if $f_\mu$ satisfies the normalization
constraint $M_3(f_\mu)=1/2\pi^2$, $f \in \mathcal{M}_{3U}$,
and we want to observe that if it happens the function $f$ is
integrable, this imply $f$ has finite
moments up to the third, i.e. it has finite
\emph{average}, \emph{variance}
and \emph{skewness} over $[0,+\infty)$.\\

It is nice to observe how the differential form
\begin{equation*}
  \frac{1}{2}\sin2\psi\ud{\psi}\ud{\phi}\ud{\theta}=
  \sin\psi\cos\psi\ud{\psi}\ud{\phi}\ud{\theta}
\end{equation*}
in the integral $\eqref{eq:corrmu}$ appears in
the definition of the SU(2) Haar Measure \cite{haarw} given
in \cite[p.~9, eq.~(4.1)]{haar} and is referenced in \cite{hopf}.\\

\section{Measures of cartesian shaped subsets of $D_h$}

We want to explore in this section some nice
properties of $SU(2)$ invariant absolutely
continuous, with respect the Lebesgue measure
$\lambda$, probability measures over $\mathcal{L}(\field{C}^2)$
for particular, cartesian product shaped, subsets of
the hypershperical domain $D_h$.\\

So, let we choose one of such absolutely continuous probability
measures $\mu$ with Radon-Nikod\'ym derivative $d\mu/d\lambda$
for which, Proposition \ref{inv3m}, there
exists a function $f_\mu \in \mathcal{M}_{3U}$,
$d\mu/d\lambda(|\vect{w}|)=p_c[\mu](|\vect{w}|)
=f_\mu(|\vect{w}|)$ a.e. over $\field{C}^2$.

Let be $\Psi$ a Lebesgue measurable subset of the real interval
$[0,\frac{\pi}{2}]$ and use it to define a cartesian product
subset of $D_h$, the measurable set
$A_{\Psi}=(0,+\infty) \times \Psi \times [0,2\pi] \times [0,2\pi]$,
i.e. in terms of cartesian coordinates, the measurable subset
of $D_c$
\begin{equation}
  \label{eq:arcval}
  T_h(A_{\Psi}) =
  \left\{
    (\vect{\alpha},\vect{\beta}) \in D_c \suchthat
    \arctan\left(\frac{|\vect{\beta}|}{|\vect{\alpha}|}\right) \in \Psi
  \right\} =
  \left[
    \arctan\left(\rho/r\right) \in \Psi
  \right]
\end{equation}
where, once again, we have taken the freedom of thinking
$\arctan(|\vect{\beta}|/|\vect{0}|)=\frac{\pi}{2}$,
when $\vect{\beta} \ne \vect{0}$, for the sake of the
notation.\\

We may apply the definition $\eqref{eq:corrmu}$ and the
Fubini's Theorem over the measure of this set to derive
\begin{equation*}\begin{split}
  \mu(T_h(A_{\Psi})) &=
  \frac{1}{2}
  \iiiint\limits_{T_h^-1(T_h(A_\psi))}
    l^3 f_\mu(l) \sin2\psi
    \ud{l}\ud{\psi}\ud{\phi}\ud{\theta}=\\
  &=
  \frac{1}{2}
  \iiiint\limits_{A_\psi}
    l^3 f_\mu(l) \sin2\psi
    \ud{l}\ud{\psi}\ud{\phi}\ud{\theta}=\\
  &=
  \frac{1}{2}
  \int_{0}^{+\infty}\int\limits_{\Psi}
  \int_{0}^{2\pi}\int_{0}^{2\pi}
    l^3 f_\mu(l) \sin2\psi
    \ud{l}\ud{\psi}\ud{\phi}\ud{\theta}=\\
  &=
  \frac{1}{2}
  \left[
    \int_{0}^{+\infty}
      l^3 f_\mu(l) \sin2\psi
  \right]
  \left[
    \int_{0}^{2\pi}\int_{0}^{2\pi}
      \ud{\phi}\ud{\theta}
  \right]
  \left[
    \int\limits_{\Psi}
      \sin2\psi
      \ud{\psi}
  \right]=\\
  &=
  \frac{1}{2}
  \Bigg[\frac{1}{2\pi^2}\Bigg]
  \Bigg[4\pi^2\Bigg]
    \left[
    \int\limits_{\Psi}
      \sin2\psi
      \ud{\psi}
  \right]=
  \int\limits_{\Psi}
    \sin2\psi
    \ud{\psi}
\end{split}\end{equation*}

It is worthing to observe we have proved the measure of
these family of subsets of the cartesian coordinates domain
$\mu(T_h(A_\Psi))=\int\limits_{\Psi}\sin2\psi\ud{\psi}$
is completely independent from the choice of the
$SU(2)$ invariant absolutely continuous probability
measure $\mu$, or, in other words, from the choice
of its Radon-Nikod\'ym derivate $f_\mu \in \mathcal{M}_{3U}$
along one radius of $\field{C}^2$.\\

This familiy of Lebesgue measurable subsets of the
cartesian coordinates domain actually is the
family of Lebesgue measurable \emph{convex cones} of $D_c$,
i.e. Lebesgue measurable subsets $K$ of $D_c$ such that
$\gamma K \subseteq K$, for each positive real number
$\gamma > 0$.\\

\begin{prop}
  \label{indip}
  The measure of the Lebesgue measurable subset
  $T_h(A_\psi)=\left[\arctan(\rho/r) \in \Psi\right]$ of $\field{C}^2$
  is
  \begin{equation*}
    \mu(T_h(A_{\Psi})) =
    \mu\left[\arctan(\rho/r) \in \Psi\right]=
    \int\limits_{\Psi}\sin2\psi\ud{\psi}
  \end{equation*}
  for each $SU(2)$ invariant, absolutely
  continuous, $\mu \ll \lambda$, probability measure
  $\mu$ over $\mathcal{L}(\field{C}^2)$ and for
  every Lebesgue measurable subset $\Psi$ of
  $[0,\frac{\pi}{2}]$.
\end{prop}

In other words Lebesgue measurable \emph{convex cones}
of $D_c$ have the same measure for each
$SU(2)$ invariant, absolutely
continuous, $\mu \ll \lambda$, probability measure
$\mu$ over $\mathcal{L}(\field{C}^2)$, the Haar
Measure of their \emph{intersection} with the unit sphere
$\left\{
(\vect{\alpha},\vect{\beta}) \in \field{C}^2
\suchthat
|\vect{\alpha}|^2 + |\vect{\beta}|^2 = 1
\right\}$
of $D_c$.\\

As a special useful case, if the set $\Psi$ is the
subinterval $[\psi_1,\psi_2]$ of $[0,\frac{\pi}{2}]$
then
\begin{equation*}
  \label{eq:measint2}
  \mu(T_h([\psi_1,\psi_2]))=
     \int_{\psi_1}^{\psi_2}
     \sin2\psi\ud{\psi}=
     \Big[\sin^2\psi\Big]_{\psi_1}^{\psi_2}=
     \sin^2\psi_2-\sin^2\psi_1
\end{equation*}
using the integral $\eqref{eq:trig1}$.\\

Will be useful
\begin{equation}
  \label{eq:measint}
  \mu(T_h([0,\psi]))=
  \mu\left[0 \le \arctan\left(\rho/r\right) \le \psi \right]=
     \sin^2\psi
  \text{,\,\,\,\,}
  0 \le \psi \le \frac{\pi}{2}
\end{equation}

Let now be $L$ a Lebesgue measurable set of $[0,+\infty)$ and
use it to define another cartesian domain of $D_h$,
the measurable set
$A_{L}=L \times [0,\frac{\pi}{2}] \times [0,2\pi] \times [0,2\pi]$,
i.e. in terms of cartesian coordinates the measurable
subset of $D_c$
\begin{equation}
  \label{eq:absval}
  T_h(A_{L}) =
  \left\{
    (\vect{\alpha},\vect{\beta}) \in D_c \suchthat
    \sqrt{|\vect{\alpha}|^2 + |\vect{\beta}|^2} \in L
  \right\} =
  \left[
     \sqrt{r^2 + \rho^2} \in L
  \right]
\end{equation}

In terms of the definition $\eqref{eq:corrmu}$ of $\mu$
applying once again the Fubini's Theorem
\begin{equation}\begin{split}
  \label{eq:indipl}
  \mu_f(T_h(A_{L})) &=
  \frac{1}{2}
  \iiiint\limits_{T_h^-1(T_h(A_L))}
    l^3 f_\mu(l) \sin2\psi \ud{l}\ud{\psi}\ud{\phi}\ud{\theta}=\\
  &=
  \frac{1}{2}
  \iiiint\limits_{A_L}
    l^3 f_\mu(l) \sin2\psi \ud{l}\ud{\psi}\ud{\phi}\ud{\theta}=\\
  &=
  \frac{1}{2}
  \int\limits_L\int_{0}^{\frac{\pi}{2}}\int_{0}^{2\pi}\int_{0}^{2\pi}
    l^3 f_\mu(l) \sin2\psi \ud{l}\ud{\psi}\ud{\phi}\ud{\theta}=\\
  &=
  \frac{1}{2}
  \left[
    \int\limits_L l^3 f_\mu(l)\ud{l}
  \right]
  \left[
    \int_{0}^{\frac{\pi}{2}} \sin2\psi \ud{\psi}
  \right]
  \left[
    \int_{0}^{2\pi}\int_{0}^{2\pi} \ud{\phi}\ud{\theta}
  \right]=\\
  &=
  \frac{1}{2}
  \Bigg[4\pi^2\Bigg]
  \Bigg[1\Bigg]
  \left[
    \int\limits_L l^3 f_\mu(l)\ud{l}
  \right]=
  2\pi^2\int\limits_L l^3 f_\mu(l)\ud{l}
\end{split}\end{equation}\\

As a special useful case, if the set L is the subinterval
$[l_1,l_2]$ of $[0,+\infty)$ then
\begin{equation*}
  \mu(T_h([l_1,l_2])) =
  2\pi^2\int\limits_{l_1}^{l_2} l^3 f_\mu(l)\ud{l}
\end{equation*}\\

It is interesting to observe the measure
of the ball $B(\vect{0},R)$ of radius $R$ of $\field{C}^2$
actually depends from the choice of the probability
measure $\mu$ and it is
\begin{equation}
  \label{eq:ball}
  \mu(B(\vect{0},R))=\mu(T_h([0,R])=
  \mu\left[ 0 \le \sqrt{r^2+\rho^2} \le R \right]=
   2\pi^2\int\limits_{0}^{R} l^3 f_\mu(l)\ud{l}
\end{equation}
for each non negative radius $R$.\\

\section{$SU(2)$ invariant continuous probability joint distributions}

Let now choose an arbitrary probability space
$(\Omega,\mathcal{F},P)$ and a cartesian coordinates domain
vector valued random variable $\vectc{W}:\Omega \to D_c$ over $\Omega$,
$\vectc{W}=(\vectc{A},\vectc{B})=(X+iY,U+iV)=(X,Y,U,V)$
with a continuous joint probability distribution defined
by an absolutely continuous normalized probability measure
$\mu:\mathcal{L}(D_c) \to [0,1]$
\begin{equation}
  P[\vectc{W} \in A] = \mu(A)
\end{equation}
for each Lebesgue measurable subset $A$ of $D_c$ and
let assume the probability distribution is invariant
for $SU(2)$ transformations
\begin{equation}
  P[\vectc{W} \in A] = P[\vectc{W} \in UA]
\end{equation}
for each unitary matrix $U \in SU(2)$ and for
every Lebesgue measurable subset $A$ of $D_c$.\\

This is equivalent to assume the absolutely continuous
probability measure $\mu$ is invariant for $SU(2)$ transformations
and, for Proposition \ref{inv3m}, its Radon-Nikod\'ym
derivative $d\mu/d\lambda$ it is equal a.e. to a
function $f_\mu \in \mathcal{M}_{3U}$ such that
\begin{equation}
  \label{eq:prob}
  P[\vectc{W} \in A] = \mu(A) =
  \frac{1}{2}
  \iiiint\limits_{T_h^-1(A)}
    l^3 f_\mu(l) \sin2\psi
    \ud{l}\ud{\psi}\ud{\phi}\ud{\theta}
\end{equation}
for every Lebesgue measurable subset $A$ of $D_c$,
using the definition $\eqref{eq:corrmu}$.\\

We may ask to ourselves which is the continuos
distribution of the real valued, non negative, random variable
$|\vectc{W}|=\sqrt{|\vectc{A}|^2+|\vectc{B}^2|}$,
absolute value of the $\field{C}^2$
vector valued random variable $\vectc{W}$.\\

To answer this question we use $\eqref{eq:indipl}$ to evaluate
\begin{equation*}\begin{split}
  P\left[|\vectc{W}| \in L\right] &=
  P\left[\sqrt{|\vectc{A}|^2 + |\vectc{B}|^2} \in L\right] =
  \mu\left[\sqrt{|\vectc{A}|^2 + |\vectc{B}|^2} \in L\right] =\\
  &=
  \mu(T_h(L))=
  2\pi^2\int\limits_L l^3 f_\mu(l)\ud{l}  
\end{split}\end{equation*}
for every Lebesgue measurable subset L of $[0,+\infty)$ and
we identify the density of the continuous probability distribution
of the real valued, non negative random variable $|\vectc{W}|$
in the real valued, non negative, integrable function
\begin{equation}
  \label{eq:density3}
  p[\mu]_{|\vectc{W}|}(l)=2\pi^2 l^3 f_\mu(l)
  \mbox{,\,\,\,\,}
  l \ge 0  
\end{equation}
in other words, using the property $\eqref{eq:ball}$
of the measure $\mu(T_h([0,l])$
\begin{equation*}
  P\left[|\vectc{W}|<l\right] =
  2\pi^2
  \int_{0}^{l} q^3 f_\mu(q) \ud{q}=
  \int_{0}^{l}p[\mu]_{|\vectc{W}|}(q)\ud{q}
  \text{,\,\,\,\,}
  l > 0
\end{equation*}
and $f_\mu \in \mathcal{M}_{3U}$ assure
the real valued, non negative, random variable $|\vectc{W}|$
has a normalized continuous probability distribution
$p[\mu]_{|\vectc{W}|}$ over $[0,+\infty)$.\\

\section{Measure of a convex cones chain}

We are interested to use this $SU(2)$
invariant, continuous probability joint distribution
$P\left[\vectc{W} \in A\right]=\mu(A)$, defined in $\eqref{eq:prob}$,
for computing the probabilities of this chain (totally ordered family)
of convex cones of the cartesian coordinates domain $D_c$, identified by a
single, real, non negative, parameter $t \ge 0$:
\begin{equation}
  \label{eq:intdef}
  A_t=\left\{{(\vect{\alpha},\vect{\beta}) \in D_c \suchthat
           0 \le |\vect{\beta}| \le t |\vect{\alpha}|} \right\}
    \text{,\,\,\,\,}
    t \in \field{R}, t \ge 0
\end{equation}
\begin{alignat*}{2}
  &t_1 \le t_2 \Leftrightarrow A_{t_1} \subset A_{t_2} \\
  &t_1 \le t_2 \Rightarrow 
  P\left[\vectc{W} \in A_{t_1}\right] \le P\left[\vectc{W} \in A_{t_2}\right]
\end{alignat*}
which defines the real valued, non negative, monotonically
non-decreasing function:
\begin{equation}
  \label{eq:intfunc}
  p(t) = P\left[\vectc{W} \in A_t\right]=
  P\left[0 \le |\vect{\beta}| \le t|\vect{\alpha}|\right]
  \text{\,\,\,\,\,}
  t \in \field{R}, t \ge 0
\end{equation}

In terms of the probability measure $\eqref{eq:prob}$ and
using our notation:
\begin{equation}
  p(t)=P\left[0 \le \rho \le tr\right]=\mu(A_t)=
  \mu(T_h([0,\arctan{t}])
  \text{\,\,\,\,\,}
  t \in \field{R}, t \ge 0
\end{equation}
and recalling $\eqref{eq:measint}$
\begin{equation}\begin{split}
  p(t)&=
  \mu_f(T_h([0,\arctan{t}])=
  \sin^2\arctan{t}=
  \left(\frac{t}{\sqrt{1+t^2}}\right)^2=\\
  &=
  \frac{t^2}{1+t^2}
  \text{\,\,\,\,\,}
  t \in \field{R}, t \ge 0
\end{split}\end{equation}

we have computed the probability measure of the chain
$\left\{A_t\right\}_{t \ge 0}$\\
\begin{equation}
  \label{eq:interes}
  p(t)=
  P\left[0 \le \rho \le tr\right]=
  \frac{t^2}{1+t^2}
  \text{\,\,\,\,\,\,\,\,\,\,}t \in \field{R},t \ge 0
\end{equation}\\
and proved $p(t)$ is independent from the choice
of the function $f_\mu \in \mathcal{M}_3$, i.e. it is
independent from the choice of the $SU(2)$
invariant, absolutely continuous probability
measure $\mu$.\\

It worths noting, recalling a previous observation,
that this independence is directly connected from the
fact the subsets of the chain $\left\{A_t\right\}_{t \ge 0}$
have a special cartesian decomposition in the
hyperspherical coordinates domain $D_h$ and from
the measure independence established in Proposition
\ref{indip}.\\

\section{The gaussian $N(0,1)$ case}

Now we choose four normal distributed, $N(0,1)$, independent,
real valued random variables $X,Y,U,V$:
\begin{equation*}
  f_X(s)=f_Y(s)=f_U(s)=f_V(s)=
    \frac{1}{\sqrt{2\pi}}
    \ex{-\frac{s^2}{2}}
\end{equation*}\\
From the definition of $\vectc{W}=(\vectc{A},\vectc{B})=(X+iY,U+iV)=
(X,Y,U,V)$ and from the independence of the real valued
random variables $X,Y,U,V$ we get the continuous joint density functions
of the vector random variables $\vectc{A}, \vectc{B}$ and $\vectc{W}$
\begin{eqnarray}
  \label{eq:casea}
  && f_{\vectc{A}}(\vect{\alpha})=f_X(x)f_Y(y)=
    \frac{1}{2\pi}\ex{-\frac{|\vect{\alpha}|^2}{2}}\\
  \label{eq:caseb}
  && f_{\vectc{B}}(\vect{\beta})=f_U(u)f_V(v)=
    \frac{1}{2\pi}\ex{-\frac{|\vect{\beta}|^2}{2}}\\
  \label{eq:case}
  && f_{\vectc{W}}(\vect{\alpha},\vect{\beta})=
    f_{\vectc{A}}(\vect{\alpha})f_{\vectc{B}}(\vect{\beta})=
    \frac{1}{4\pi^2}\ex{-\frac{|\vect{\alpha}|^2+|\vect{\beta}|^2}{2}}
\end{eqnarray}

The joint density function $\eqref{eq:case}$ over the
cartesian coordinates domain $D_c$ is the Radon-Nikod\'ym derivative
$d\mu/d\lambda$ of the absolutely continuous measure $\mu$ over
$\mathcal{L}(D_c)$ defined by
\begin{equation}
  \mu(A) =
  \int\limits_{A}f_{\vectc{W}}(\vect{w})\vdm{4}{\lambda}{w}=
  \int\limits_{A}
    \frac{1}{4\pi^2}
    \ex{-\frac{|\vect{w}|^2}{2}}
    \vdm{4}{\lambda}{w}
  \mbox{,\,\,\,\,}A \in \mathcal{L}(D_c)
\end{equation}
and the real valued, non negative, integrable function
$f_\mu:[0,+\infty) \to [0,+\infty)$
\begin{equation}
  f_\mu(l)=\frac{1}{4\pi^2}
  \ex{-\frac{l^2}{2}}
  \mbox{,\,\,\,\,} l \ge 0
\end{equation}
has exactly the property we expressed in definition 
$\eqref{eq:corrf}$
\begin{equation}
  \label{eq:corrfg}
  f_\mu(l)=\frac{d\mu}{d\lambda}(l,\vect{0})=
  f_{\vectc{W}}(l,\vect{0})=
   \frac{1}{4\pi^2}\ex{-\frac{l^2}{2}}
  \mbox{,\,\,\,\,} l \ge 0
\end{equation}
and what is left, for Proposition \ref{inv3m},
to verify the absolutely continuous
measure $\mu$ is an SU(2) transformations invariant \emph{finite}
measure, is to check $f_\mu$ has finite third moment $M_3(f_\mu)$.

Using the integral $\eqref{eq:exp2}$
\begin{equation*}\begin{split}
  M_3(f_\mu) &=
  \int_{0}^{+\infty}
    l^3 f_\mu(l)\ud{l}=
  \frac{1}{4\pi^2}
  \int_{0}^{+\infty}
    l^3\ex{-\frac{l^2}{2}}\ud{l}=
  \frac{1}{4\pi^2}
  \int_{0}^{+\infty}
    \frac{1}{2}l\ex{-\frac{l}{2}}\ud{l}=\\
  &=
  \frac{1}{4\pi^2}\frac{1}{2}[4] = \frac{1}{2\pi^2} < +\infty
\end{split}\end{equation*}
it is easy to check the function $f_\mu$ is actually a
member of the subclass $\mathcal{M}_{3U}$ because it
satisfies the \emph{normalization constraint} $M_3(f_\mu)=1/2\pi^2$,
how it has to expected from the fact the density
function $f_{\vectc{W}}$ is already normalized, being
a product of normalized gaussians.

Using Proposition \ref{inv3m} we are now legimate to
conclude the measure $\mu$ is an absolutely continuos
probability measure, $SU(2)$ invariant probability
joint distribution of the random variable $\vectc{W}$
and we may apply the definition $\eqref{eq:corrmu}$
for an explicit definition of the measure in
hyperspherical coordinates
\begin{equation*}
  \mu(A)=
  \frac{1}{2}
  \iiiint\limits_{T_h^-1(A)}
    l^3 \ex{-\frac{l^2}{2}} \sin2\psi
    \ud{l}\ud{\psi}\ud{\phi}\ud{\theta}
\end{equation*}  
for each Lebesgue measurable subset $A$ of $D_c$.

This verification allows us to use $\eqref{eq:density3}$
to compute the density function $f_{|\vectc{W}|}$ of the
real valued, non negative, random variable $|\vectc{W}|$
\begin{equation}
  \label{eq:absw}
  f_{|\vectc{W}|}(l) =
  2\pi^2 l^3 f_\mu(l) =
  2\pi^2 \frac{1}{4\pi^2} l^3 \ex{-\frac{l^2}{2}}=
  \frac{1}{2} l^3 \ex{-\frac{l^2}{2}}
  \text{,\,\,\,\,}
  l \ge 0
\end{equation}
and to immediately state that the result
achieved in $\eqref{eq:interes}$ for the chain
$\left\{A_t\right\}_{t \ge 0}$ defined in 
$\eqref{eq:intdef}$ holds for this gaussian continuous
probability distribution.\\

But, as a further verification, we want to follow
a different path to prove $\eqref{eq:interes}$ and $\eqref{eq:absw}$
hold for such a distribution.\\

From the joint density functions $\eqref{eq:casea}$ and $\eqref{eq:caseb}$
of the vector random variables $\vectc{A}=(X,Y)$ and $\vectc{B}=(U,V)$
we know the real, non negative, valued random variables of their
absolute values $|\vectc{A}|$ and $|\vectc{B}|$ are Rayleigh(1)
distributed \cite{rayleigh}, i.e.:
\begin{equation*}
  f_{|\vectc{A}|}(s) =
  f_{|\vectc{B}|}(s) =
  s\ex{-\frac{s^2}{2}}
  \text{,\,\,\,\,}
  s \ge 0\\ 
\end{equation*}

Now we choose a measurable subset $S$ of $\field{R}^2$
of couples $(r,\rho)$, such that $r \ge 0, \rho \ge 0$
$r > 0$ or $\rho >0$, and evaluate the probability the couple
$(|\vectc{A}|,|\vectc{B}|)$ of real, non negative, numbers
belongs to $S$.

Using our assumptions, the double polar
transformation $T_h\text{\,\,}\eqref{eq:dpolar}$, the integral
evaluation rule $\eqref{eq:ruled}$ and the Fubini's Theorem
\begin{multline*}
  P\left[
    (|(\vectc{A}|,|\vectc{B}|) \in S
  \right]=
  \iiiint\limits_{\left[(|(\vectc{A}|,|\vectc{B}|) \in S\right]}
    f_{\vectc{W}}(\vect{\alpha},\vect{\beta})
    \vd{2}{\alpha}\vd{2}{\beta}=\\
  =
  \iint\limits_{S}\int_{0}^{2\pi}\int_{0}^{2\pi}
    r\rho\frac{1}{4\pi^2}\ex{-\frac{r^2+\rho^2}{2}}
    \ud{r}\ud{\rho}\ud{\phi}\ud{\theta}=
  \iint\limits_{S}
    r\rho\ex{-\frac{r^2+\rho^2}{2}}
    \ud{r}{\rho}=\\
  =
  \iint\limits_{S}
    \Big[
      r\ex{-\frac{r^2}{2}}
    \Big]
    \Big[
      \rho\ex{-\frac{\rho^2}{2}}
    \Big]
    \ud{r}\ud{\rho}=
  \iint\limits_{S}
    f_{|\vectc{A}|}(r)f_{|\vectc{B}|}(\rho)
    \ud{r}\ud{\rho}
\end{multline*}
for every of such measurable subsets $S$ of $\field{C}$,
in other words, in terms of density functions:
$f_{(|\vectc{A}|,|\vectc{B}|)}(r,\rho)=
f_{|\vectc{A}|}(r)f_{|\vectc{B}|}(\rho)$ a.e. for
$r \ge 0, \rho \ge 0$.\\

We have just proved the real valued non negative random
variables $|\vectc{A}|,|\vectc{B}|$ are independent
random variables and we may now use a standard polar coordinates
$(l,\psi)$ transformation to evaluate
\begin{multline*}
  P\left[
    \sqrt{|\vectc{A}|^2+|\vectc{B}|^2} \in L
  \right]=
  \iint\limits_{\left[\sqrt{|\vectc{A}|^2+|\vectc{B}|^2} \in L\right]}
    f_{|\vectc{A}|}(r)f_{|\vectc{B}|}(\rho)
    \ud{r}\ud{\rho}=\\
  =
  \iint\limits_{\left[\sqrt{|\vectc{A}|^2+|\vectc{B}|^2} \in L\right]}
    r\rho\ex{-\frac{r^2+\rho^2}{2}}
    \ud{r}\ud{\rho}=
  \int\limits_{L}\int_{0}^{\frac{\pi}{2}}
    (l\cos\psi)(l\sin\psi)
    \left(\ex{-\frac{l^2}{2}}\right)
    l\ud{l}\ud{\psi}=\\
  =
  \int\limits_{L}\int_{0}^{\frac{\pi}{2}}
    \frac{1}{2}l^3\ex{-\frac{l^2}{2}}\sin2\psi
    \ud{l}\ud{\psi}=
  \Bigg[
    \int_{0}^{\frac{\pi}{2}}
      \sin2\psi
      \ud{\psi}
  \Bigg]
  \Bigg[
    \int\limits_{L}
      \frac{1}{2}l^3\ex{-\frac{l^2}{2}}
      \ud{l}
  \Bigg]=\\
  =
  \int\limits_{L}
    \frac{1}{2}l^3\ex{-\frac{l^2}{2}}
    \ud{l}
\end{multline*}
for every measurable subset $L$ of $(0,+\infty)$. 
We have just recovered the same density function
obtained in $\eqref{eq:absw}$.\\

From the definition $\eqref{eq:case}$
of $f_{\vectc{W}}(\vect{\alpha},\vect{\beta})$,
using our double polar coordinates transformation
$T_h\text{\,\,}\eqref{eq:dpolar}$, integral evalutation
rule $\eqref{eq:ruled}$ and the Fubini's Theorem, we may directly
compute the probabilities of the subsets $\left\{A_t\right\}_{t \ge 0}$
\begin{equation*}\begin{split}
  p(t)&=P[W \in A_t]=P[0 \le |\vect{\beta}| \le t|\vect{\alpha}|]=
  \iiiint\limits_{[0 \le |\vect{\beta}| \le t\|\vect{\alpha}|]}
    f_{\vectc{W}}(\vect{\alpha},\vect{\beta})
    \vd{2}{\alpha}\vd{2}{\beta}=\\
  &=
  \int\limits\limits_{[0 \le \rho \le tr]}
  \int_{0}^{2\pi}\int_{0}^{2\pi}
    r\rho \frac{1}{4\pi^2}\ex{-\frac{r^2+\rho^2}{2}}
    \ud{r}\ud{\rho}\ud{\phi}\ud{\theta}=\\
  &=
  \int_{0}^{+\infty}\int_{0}^{tr}
    r\rho \ex{-\frac{r^2+\rho^2}{2}}
    \ud{\rho}\ud{r}=
  \int_{0}^{+\infty}
  \Bigg[
    r\ex{-\frac{r^2}{2}}
    \int_{0}^{tr}
    \rho\ex{-\frac{\rho^2}{2}}
    \ud{\rho}
  \Bigg]
  \ud{r}=\\
  &=
  \int_{0}^{+\infty}
    r\ex{-\frac{r^2}{2}}
    \Big[
      -\ex{-\rho}
    \Big]_{0}^{\frac{t^2r^2}{2}}
    \ud{r}=
  \int_{0}^{+\infty}
    r\ex{-\frac{r^2}{2}}
    \left(
      -\ex{\frac{t^2r^2}{2}}+1
    \right)
    \ud{r}=\\
  &=
    \int_{0}^{+\infty}
    r\ex{-\frac{r^2}{2}}
    \ud{r}
    -
    \int_{0}^{+\infty}
      r\ex{-\frac{r^2}{2}}
      \ex{\frac{t^2r^2}{2}}
      \ud{r}=
    1
    -
    \int_{0}^{+\infty}
      r\ex{-(1+t^2)\frac{r^2}{2}}
      \ud{r}=\\
  &=
    1
    -
    \int_{0}^{+\infty}
      \frac{1}{1+t^2}\ex{-r}
      \ud{r}=
    1
    -
    \frac{1}{1+t^2}
    \Big[
      -\ex{-r}
    \Big]_{0}^{+\infty}=
    1-\frac{1}{1+t^2}=
    \frac{t^2}{1+t^2}
 \end{split}\end{equation*}
for every $t \ge 0$ and
we have just proved again the property $\eqref{eq:interes}$
of the function $p(t)$ holds for the gaussian $N(0,1)$ case.\\

\section{Conclusion}

Once we know the function $p(t)$, defined in $\eqref{eq:intfunc}$,
is independent from the $SU(2)$ transformations invariant
probability continuous joint distribution $\mu_f$,
for every $f \in \mathcal{M}_{3U}$, it would have been completely
justified to evaluate $p(t)$ by picking the gaussian shaped
probability continuous joint distribution $\eqref{eq:case}$,
as done in \cite[p.~10, eq.~(65)/(66)]{quantum}.

But we may directly achieve the same result using the
general property we got in $\eqref{eq:interes}$
\begin{equation*}\begin{split}
  p\left(|\vect{a}|/|\vect{b}|\right) &=
  P\Bigg[ 0 \le \rho \le \frac{|\vect{a}|}{|\vect{b}|}r\Bigg]=
  P\Big[ |\vect{a}|r \ge |\vect{b}|\rho \Big]=
  P\Big[ |\vect{a}||\vectc{A}| \ge |\vect{b}||\vectc{B}| \Big]=\\
  &=
  \frac{
    \frac{|\vect{a}|^2}{|\vect{b}|^2}
  }
  {
    1 + \frac{|\vect{a}|^2}{|\vect{b}|^2}
  }=
  \frac{
    |\vect{a}|^2
  }
  {
    |\vect{a}|^2 + |\vect{b}|^2
  }
\end{split}\end{equation*}
for each $\vect{a}, \vect{b} \in \field{C}^2, \vect{b} \ne \vect{0}$,
independently from the choice of the $SU(2)$ transformations
invariant continuous probability joint distribution of
the random variable $\vectc{W}=(\vectc{A},\vectc{B})$.

This is essentially proved in Propositions $\ref{inv2}$ and
$\ref{inv3m}$ using the existence of the Haar probability
measure over compact Hausdorff topological groups, a rather
deep result of Measure Theory.


\newpage

\begin{thebibliography}{1}
\bibitem{aocc}
  \emph{Axiom of countable choice}\\
  Wiki page about the Axiom of countable choice,\\
  \url{http://en.wikipedia.org/wiki/Axiom_of_countable_choice}

\bibitem{haarw}
  \emph{Haar Measure}\\
  Wiki page about the Haar Measure,\\
  \url{http://en.wikipedia.org/wiki/Haar_measure}

\bibitem{hopf}
  \emph{Hopf coordinates}\\
  Wiki page about the 3-sphere,\\
  \url{http://en.wikipedia.org/wiki/3-sphere#Hopf_coordinates}

\bibitem{moment}
  \emph{Moment}\\
  Wiki page about the Moment (mathematics),\\
  \url{http://en.wikipedia.org/wiki/Moment_(mathematics)}

\bibitem{rayleigh}
  \emph{Rayleigh distribution}\\
  Wiki page about the Rayleigh distribution,\\
  \url{http://en.wikipedia.org/wiki/Rayleigh_distribution}

\bibitem{bartle2}
  \emph{The Elements of Integration and Lebesgue Measure}\\
  Robert G. Bartle\\
  John Wiley \& Sons, 1995

\bibitem{fremlin2}
  \emph{Measure Theory, Volume 2}\\
  D.H. Fremlin\\
  \url{http://www.essex.ac.uk/maths/people/fremlin/mt.htm}
  
\bibitem{fremlin4i}
  \emph{Measure Theory, Volume 4, Part I}\\
  D.H. Fremlin\\
  \url{http://www.essex.ac.uk/maths/people/fremlin/mt.htm}
  
\bibitem{haar}
  \emph{Representations of SU(2) and Jacobi polynomials}\\
  Tom H. Koornwinder, thk@science.uva.nl,\\
  \url{http://staff.science.uva.nl/~thk/edu/orthopoly.pdf}

\bibitem{naimark}
  \emph{Normed Rings}\\
  M.A. Naimark\\
  Wolters-Noordhoff, 1970

\bibitem{rudin}
  \emph{Real and Complex Analysis}\\
  Walter Rudin,\\
  McGraw-Hill, 1970
  
\bibitem{quantum}
  \emph{A Theory of Quantum Observation and the Emergence of the Born Rule}\\
  Andreas O. Tell,\\
  \url{http://arxiv.org/abs/1205.0293}

\end{thebibliography}
\end{document}